\documentclass[preprint,3p,11pt]{elsarticle}
\usepackage{setspace}
\usepackage{booktabs}
\usepackage{graphicx}
\usepackage[short,c2,nocomma]{optidef}
\usepackage{xcolor}
\usepackage{amsmath, amsthm, amsfonts, amssymb}
\usepackage{enumitem}

\usepackage[colorlinks=true,linkcolor=blue, citecolor=blue, urlcolor=blue, breaklinks]{hyperref}

\newtheorem*{remark}{Remark}

\newlength{\digitwidth}
\begin{document}

\begin{frontmatter}

\title{ICNN-enhanced 2SP: Leveraging input convex neural networks for solving two-stage stochastic programming}

\author[1]{Yu Liu} 
\author[2,1]{Fabricio Oliveira\corref{cor}}
\cortext[cor]{Corresponding author.}
\ead{fabol@dtu.dk}
\author[3]{Jan Kronqvist}

\affiliation[1]{organization={Department of Mathematics and Systems Analysis, School of Science, Aalto University}, city={Espoo}, country={Finland}}
\affiliation[2]{organization={DTU Management, Technical University of Denmark}, city={Lyngby}, country={Denmark}}
\affiliation[3]{organization={Department of Mathematics, KTH Royal Institute of Technology}, city={Stockholm}, country={Sweden}}

\begin{abstract}
Two-stage stochastic programming (2SP) offers a basic framework for modelling decision-making under uncertainty, yet scalability remains a challenge due to the computational complexity of recourse function evaluation. 
Existing learning-based methods like Neural Two-Stage Stochastic Programming (Neur2SP) employ neural networks (NNs) as recourse function surrogates but rely on computationally intensive mixed-integer programming (MIP) formulations.
We propose ICNN-enhanced 2SP, a method that leverages Input Convex Neural Networks (ICNNs) to exploit linear programming (LP) representability in convex 2SP problems.
By architecturally enforcing convexity and enabling exact inference through LP, our approach eliminates the need for integer variables inherent in the conventional MIP-based formulation while retaining an exact embedding of the ICNN surrogate within the 2SP framework. This results in a more computationally efficient alternative, and we show that good solution quality can be maintained.
Comprehensive experiments reveal that ICNNs incur only marginally longer training times while achieving validation accuracy on par with their standard NN counterparts. Across benchmark problems, ICNN-enhanced 2SP often exhibits considerably faster solution times than the MIP-based formulations while preserving solution quality, with these advantages becoming significantly more pronounced as problem scale increases. For the most challenging instances, the method achieves speedups of up to 100$\times$ with solution quality superior to MIP-based formulations.
\end{abstract}



\end{frontmatter}

\section{Introduction}

Two-stage stochastic programming (2SP) is a cornerstone modelling framework for decision-making under uncertainty, enabling optimal resource allocation in the face of unpredictable future events. In this, initial (first-stage) decisions are made before uncertainty is resolved, followed by adjustment (recourse, or second-stage) decisions once the random variables have unfolded. Industries such as power~\citep{romisch_recent_2010} and petrochemical sectors~\citep{khor_two-stage_2008} rely on 2SP to incorporate uncertainty in their optimisation modelling, utilising it as a standard tool in planning activities. 

Sample Average Approximation (SAA)~\citep{shapiro_simulation-based_1998} is a foundational technique for 2SP, in which the expected value of the second-stage problem is replaced with a finite-scenario approximation. Though algorithmically straightforward to implement, SAA faces a fundamental trade-off: insufficient scenarios may fail to accurately represent the random phenomena, while larger samples increase computational demands. 
To alleviate this limitation, scenario reduction techniques~\citep{heitsch_scenario_2003} are often applied, but they risk omitting critical scenarios, undermining solution robustness. 
These shortcomings have prompted the adoption of specialised decomposition approaches, including Benders-based~\citep{sherali_modification_2002} and Lagrangian decomposition-based methods~\citep{caroe_dual_1999}, which decompose the problem into smaller, more manageable subproblems. 
By exploiting problem structure (e.g., separability and duality), these methods mitigate SAA’s scalability challenges and reduce dependency on scenario sampling~\citep{kucukyavuz_introduction_2017}. However, their efficacy remains contingent on problem-specific properties, such as linearity or convexity, and often requires meticulous algorithmic tuning. 
Despite developments in decomposition methods~\citep{varga_speeding_2024}, these limitations persist, driving exploration of alternative approaches that leverage machine learning (ML).

Recent advances in learning-based methods for 2SP have explored surrogate models to approximate complex value functions. 
Several studies have proposed neural network (NN)-based approaches for learning linear programming (LP) value functions, aiming to efficiently solve two- and multi-stage stochastic programming problems~\citep{dai_neural_2021, bae_deep_2023, larsen_fast_2024, anh-nguyen_learning_2024}.
Other works have focused on learning-based scenario reduction to create smaller representative scenario sets~\citep{bengio_learning-based_2020, wu_learning_2021, wu_hgcn2sp_2024}. 
Another line of research decomposes 2SP into multi-agent structures trained on sampled scenarios via deep reinforcement learning, with the training phase emerging as the primary computational bottleneck~\citep{hubbs_or-gym_2020, yilmaz_deep_2024, stranieri_combining_2024}. 
While these approaches demonstrate the advantages of surrogate modelling, they are relatively less generalisable, constrained by specific problem components or reliance on bespoke architectures. In contrast, Neur2SP (Neural Two-Stage Stochastic Programming)~\citep{dumouchelle_neur2sp_2022} presents a general framework by employing NNs to approximate the expected recourse function, which is motivated by the embedding of a trained NN surrogate~\citep{fischetti_deep_2018}.  
\citet{kronqvist_alternating_2024} further improved Neur2SP through dynamic sampling and retraining to refine the surrogate in regions of interest, a principle also embodied in the adaptive surrogate optimisation framework of~\citet{liu_simulator-based_2025}.
However, Neur2SP’s generality comes at a cost: embedding NNs with rectified linear unit (ReLU) activations into optimisation models requires mixed-integer programming (MIP) formulations that scale with the network size, leading to computational intractability for deeper or wider architectures. 
This tension between generality and efficiency underscores the need for architectures that balance expressivity with tractability.

Input Convex Neural Networks (ICNNs)~\citep{amos_input_2017} have emerged as a promising alternative for surrogate modelling due to their ability to fit general convex functions while retaining the flexibility and expressivity of NNs. Unlike typical feed-forward networks, ICNNs are architecturally constrained to ensure convexity with respect to their inputs, enabling exact inference through LP formulations~\citep{duchesne_supervised_2021}. This property renders them particularly advantageous for embedding within large-scale optimisation frameworks.
Although ICNNs cannot model non-convex recourse functions, a limitation compared to ReLU NNs, many real-world 2SP problems in inventory management, energy dispatch, and portfolio optimisation inherently exhibit convex structure~\citep{agrawal_learning_2019, zhong_real-time_2021, luxenberg_portfolio_2024}. 
Despite their promise, prior research has predominantly focused on applying ICNNs for function approximation and system modelling, such as learning Optimal Power Flow (OPF) value functions~\citep{chen_data-driven_2020, rosemberg_learning_2024} or representing system dynamics in optimal control~\citep{yang_optimization-based_2021, bunning_input_2021, jiang_path-following_2022, wang_fast_2024}. 
However, their potential in surrogate-based optimisation, specifically for approximating convex recourse functions in 2SP, remains underexplored. In particular, convex neural architectures have been studied for modelling LP value functions in two-stage MIPs with continuous recourse by~\citet{anh-nguyen_learning_2024}. However, they rely on a specialised neural architecture tailored to structural properties of LPs and on a compact heuristic reformulation, rather than on general, easily deployable ICNNs. 
As such, our work bridges this gap by unlocking efficient and scalable solutions for 2SP using ICNN-based surrogates, leveraging the inherent convexity and computational tractability of LP.

\begin{figure}[b]
    \centering
    \includegraphics[width=0.70\textwidth]{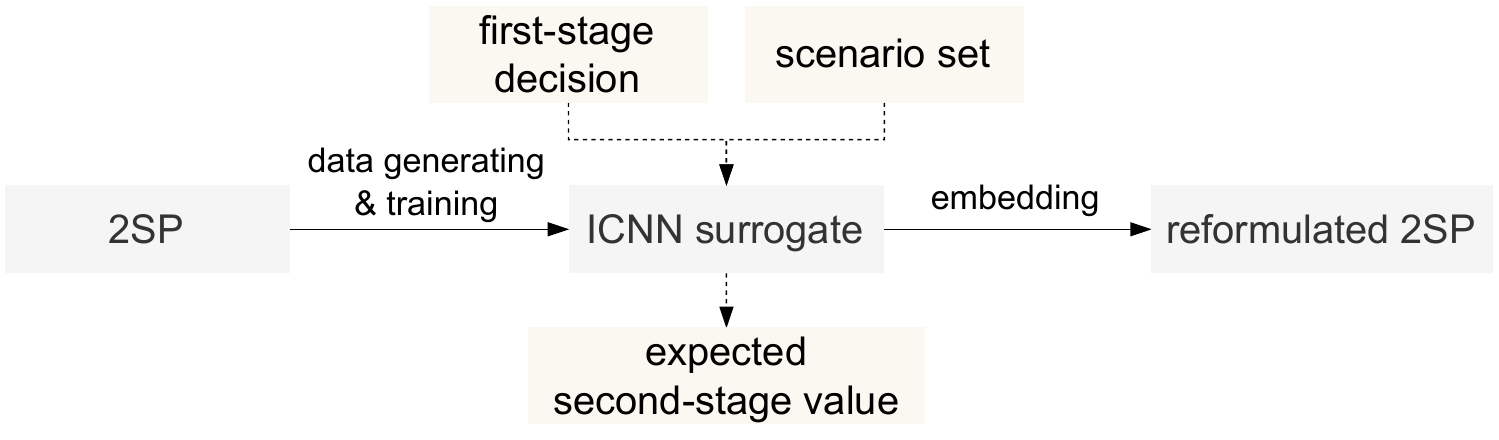}
    \caption{Overview of the proposed optimisation methodology.}
    \vspace{-6pt}
    \label{fig:intro}
\end{figure}

In particular, we present several interconnected contributions within a surrogate-based optimisation framework by integrating ICNNs into the 2SP paradigm, as we illustrate in Figure~\ref{fig:intro}. 
First, we exploit the LP representability of ICNNs to eliminate integer variables inherent to conventional MIP-based embeddings. 
Second, we propose the first integration of ICNNs into 2SP, training them to approximate the expected value function while preserving convexity.  
Third, extensive numerical experiments show that our method achieves solution quality comparable to Neur2SP while reducing computational time, especially for large-scale instances, thereby enabling application to time-critical decision problems.

The paper proceeds as follows: 
Section~\ref{sec2} reviews 2SP and ICNN fundamentals; 
Section~\ref{sec3} analyses the convexity of the surrogate target;
Section~\ref{sec4} details the process of ICNN surrogate training and embeddings; 
Section~\ref{sec5} presents the numerical experiments involving three different 2SP problems; 
Section~\ref{sec6} concludes with key findings and opening research directions.

\section{Preliminaries}\label{sec2}

This section introduces the general formulation of 2SP and the properties of ICNNs, which form the foundation of our method.

\subsection{Two-stage stochastic programming}\label{sec2_1}

Let $x \in \mathbb{R}^n$ denote first-stage decisions, $y \in \mathbb{R}^{n_y}$ the recourse (second-stage) decisions, and $\xi$ the relevant random variable. The stochastic components of the second-stage data $q(\xi) \in \mathbb{R}^{n_y}$, $W(\xi) \in \mathbb{R}^{m_y \times n_y}$, $h(\xi) \in \mathbb{R}^{m_y}$ and $T(\xi) \in \mathbb{R}^{m_y \times n}$ are determined by the realisation of $\xi$. The general linear 2SP formulation is%
\begin{mini!}
    {x}{c^\intercal x + \mathcal{Q}(x)}{\label{eq:init_1}}{}
    \addConstraint{Ax}{= b,\label{eq:X_1}}{}
    \addConstraint{x}{\geq 0,\label{eq:X_2}}{}   
\end{mini!}
where $A \in \mathbb{R}^{m \times n}$, $b \in \mathbb{R}^m$, $c \in \mathbb{R}^n$, $\mathcal{Q}(x)=\mathbb{E}_\xi [Q(x,\xi)]$ and the second-stage problem is
\begin{equation}
    Q(x,\xi)=\min_y \left\{   
    q(\xi)^\intercal y\mid W(\xi)y=h(\xi)-T(\xi)x,\ y\geq 0
    \right\},
    \label{eq:init_2}
\end{equation} 
which evaluates to the optimal second-stage cost under realisation $\xi$ given the first-stage decision $x$. We assume relatively complete recourse, ensuring $Q$ is feasible for any $x$ satisfying \eqref{eq:X_1} and \eqref{eq:X_2}.

The decision process follows $x \rightarrow \xi \rightarrow y(x,\xi)$, where $y$ is chosen to minimise $Q(x,\xi)$. We assume $\xi$ has a known probability distribution with support $\Xi$. In practical implementations, $\Xi$ is treated as a discrete and finite set where each realisation $\xi_s \in \Xi$ for $s \in S \equiv \{1, \ldots, \lvert\Xi\rvert \}$ represents a scenario.
This yields the deterministic equivalent formulation with a finite number of variables and constraints:%
\begin{mini!}
  {x}{c^\intercal x + \sum_{s\in S} p_s q_s^\intercal y_s}{\label{eq:ef}}{}
  \addConstraint{Ax }{= b,}{x \geq 0}
  \addConstraint{T_s x + W_s y_s }{= h_s,\quad}{\forall s \in S}
  \addConstraint{y_s }{\geq 0,}{\forall s \in S,}
\end{mini!}
where $p_s$ is the probability associated with scenario $s$. 
This so-called \textit{extensive form} (EF)~\eqref{eq:ef} provides a tractable formulation. However, it suffers from scalability issues as both the number of variables and constraints increase linearly with the cardinality of the scenario set $\lvert\Xi\rvert$. As such, for variants involving mixed-integer decision variables, the computational demand increases substantially for larger scenario sets.  

\subsection{ICNN architectures}\label{sec2_2}

\begin{figure}[t!]
    \centering
    \includegraphics[width=0.47\textwidth]{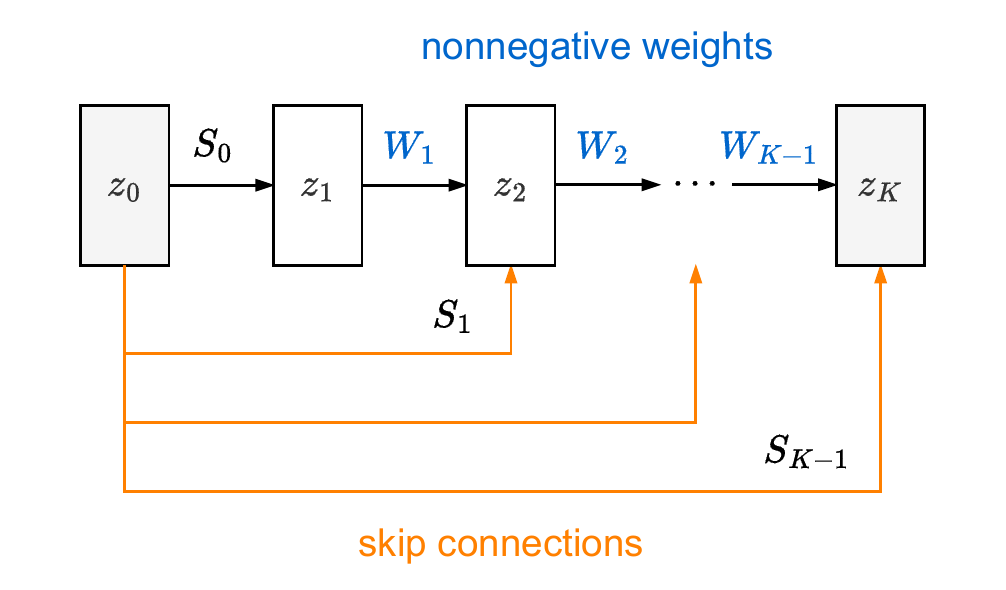}
    \caption{Architecture of a fully connected ICNN.}
    \vspace{-6pt}
    \label{fig:icnn_arch}
\end{figure}

We consider a fully connected $K$-layer ICNN, as illustrated in Figure~\ref{fig:icnn_arch}, to construct a scalar-valued NN $\hat f(z_0;\theta)$ that approximates a target function $f$. Here, $z_0$ denotes an input vector and $\theta$ represents network parameters. 
The activation for each layer can be expressed as:
\begin{equation}
    z_{k+1} = \sigma \left( W_k z_k + S_k z_0 + b_k\right), \ k = {0, \dots, K-1},
\label{eq:fcicnn}
\end{equation}
where $z_k$ for $k \geq 1$ is obtained by applying an activation function $\sigma(\cdot)$, with initial condition $W_0 \equiv 0$. The output of the network is $\hat f(z_0;\theta) = z_K$ with $\theta =\{W_{1:K-1}, S_{0:K-1}, b_{0:K-1}\}$.

The surrogate function $\hat f$ is guaranteed to be convex in $z_0$ provided that $W_{1:K-1}$ are non-negative and the activation function $\sigma(\cdot)$ (ReLU in our case) is convex and non-decreasing~\citep{amos_input_2017}. This follows standard results relating to convexity-preserving operations applied to convex functions: (i) non-negative sums of convex functions are convex; (ii) the composition of a convex function with a convex non-decreasing function retains convexity.
Another relevant feature of ICNNs is the presence of \textit{skip connections} parameterised by $S_k$, which directly propagates the input $z_0$ to deeper layers. This design aims to mitigate the restrictive effects of the non-negativity constraints imposed on $W_k$.
Notably, these skip connections preserve convexity as the sum of a convex function (from previous layers) and a linear function (from the skip connections) remains convex.

\subsection{Exact inference in ICNNs}\label{sec2_3}

Unlike standard NNs, which generate predictions solely through feedforward propagation of inputs across layers, the original formulation of ICNNs requires the solution of an optimisation problem as a subroutine to derive predictions.
For a fully connected ICNN as described in~\eqref{eq:fcicnn}, inference for a given input $z_0$ is mathematically equivalent to solving the following LP problem~\citep{amos_input_2017}:%
\begin{mini!}
    {z_k}{z_K}{\label{eq:icnn_lp}}{}
    \addConstraint{z_{k+1} }{\geq W_k z_k + S_k z_0 + b_k,\quad \label{eq:icnn_c1}}{k=0, \ldots ,K-1}
    \addConstraint{z_k }{\geq 0,\quad \label{eq:icnn_c2}}{k=1, \ldots ,K-1.}
\end{mini!}

In this formulation, the non-linearity of ReLU activations and layer-wise equality constraints are transformed into inequalities via an \textit{epigraph} representation, exploiting the fact that ReLU is a piecewise-linear function. One may note that the minimisation objective, together with the convex inequality constraints, enforces that, for each neuron, at least one of the constraints in~\eqref{eq:icnn_c1} or~\eqref{eq:icnn_c2} becomes active at the optimal solution, recovering the original network’s feedforward behaviour.
This transformation provides a principled framework for embedding a trained ICNN surrogate into larger optimisation problems, as it offers not only an exact representation of the NN but also preserves computational tractability inherent in LP.

\section{Convexity analysis}\label{sec3}

This section establishes the theoretical foundation for using ICNNs as surrogates for the second-stage cost $\mathcal{Q}(x)$ in 2SP. 
Since ICNNs inherently preserve convexity, their effectiveness hinges on whether $\mathcal{Q}(x)$ is convex, or whether it can be, to a sufficient degree, approximated by a piecewise-linear convex function.  
We characterise conditions under which convexity of $\mathcal{Q}(x)$ can be guaranteed, and show that for certain nonconvex second-stage problems with specific structure, $\mathcal{Q}(x)$ exhibits quasi-convexity. 
This motivates ICNNs as a strong class of approximators of $\mathcal{Q}(x)$ across a range of practically relevant settings.

\subsection{Case 1: continuous first- and second-stage variables}
When $x \in \mathbb{R}^n_+$ and $y_s \in \mathbb{R}^{n_y}_+$, $\forall s \in S$, the second-stage problem $Q(x, \xi_s)$ is an LP. By LP duality, we have that
\begin{equation}
    Q(x,\xi_s)=\max_\pi \left\{   
   \pi^\intercal(h_s-T_sx) \mid W^\intercal_s\pi \leq q_s
    \right\},
    \label{eq:init_dl1}
\end{equation} 
where $\pi$ represents the dual variable associated with the equality constraint $W_sy = h_s - T_sx$. The dual objective is affine in $x$ and the dual feasible set
\begin{equation}
D_s(x) = \left\{\pi \in \mathbb{R}^{m_y} \mid W^\intercal_s\pi \leq q_s \right\}
\end{equation}
constitutes a convex polyhedron. By the fundamental theorem of LP, the maximum of the dual objective is attained at an extreme point (vertex) of $D_s(x)$, given that the relatively complete recourse assumption guarantees that \eqref{eq:init_dl1} is feasible and bounded.
Let $\pi_1, \pi_2, \ldots, \pi_K$ denote the finite set of extreme points of $D_s(x)$. Thus:
\begin{equation}
    Q(x,\xi_s)=\max_{k\in\{1,\ldots,K\}} \left\{   
   -\pi^\intercal_k T_sx +\pi^\intercal_kh_s
    \right\}.
    \label{eq:init_dl2}
\end{equation} 
Each term $-\pi^\intercal_k T_sx +\pi^\intercal_kh_s$ is affine in $x$, so $Q(x,\xi_s)$ is the pointwise maximum of affine functions. As affine functions are convex and the pointwise maximum preserves convexity~\citep[Section 3.2.3]{boyd_convex_2004}, $Q(x, \xi_s)$ is convex in $x$. The expected second-stage cost $\mathbb{E}_{\xi}[Q(x,\xi)] = \sum_{s\in S} p_s Q(x, \xi_s)$ is a convex combination of convex functions $Q(x, \xi_s)$, since $0 \le p_s \le 1$, $\forall s \in S$, and $\sum_{s \in S}p_s = 1$. Thus, $\mathcal{Q}(x)$ is convex in $x$.

\begin{remark}[Extension to nonlinear convex recourse]
    The convexity of $\mathcal{Q}(x)$ extends to second-stage problems with nonlinear convex constraints, provided that appropriate constraint qualification (CQ) conditions are satisfied, such as Slater's CQ condition~\citep[Section 5.2.3]{boyd_convex_2004}, since strong duality still applies under CQ. 
    In this more general setting, the dual formulation transforms $Q(x,\xi_s)$ into the pointwise maximum of convex functions in $x$. Following similar duality argument, $\mathcal{Q}(x)$ preserves its convexity.
    Furthermore, if the first-stage problem includes a nonlinear convex objective, along with nonlinear convex inequality constraints and affine equality constraints, the recourse function remains convex in $x$~\citep[Section 3.5]{birge_introduction_2011}.
\end{remark}

\subsection{Case 2: mixed-integer first- and continuous second-stage variables}
Consider $x \in \mathcal{X}$ where $\mathcal{X} = \{ x \in \mathbb{Z}^n_+ \mid Ax=b \}$ and $y_s \in \mathbb{R}^{n_y}_+$. 
Let us define $\mathcal{X}_\text{relax} = \{ x \in \mathbb{R}^n_+ \mid Ax=b \}$ as the continuous relaxation of $\mathcal{X}$, and $\operatorname{conv}(\mathcal{X})$ as the convex hull of $\mathcal{X}$.
From Case 1, we established that $\mathcal{Q}(x)$ is convex over $\mathcal{X}_\text{relax}$. 
Since $\operatorname{conv}(\mathcal{X}) \subset \mathcal{X}_\text{relax}$, the function $\mathcal{Q}(x)$ is also convex over $\operatorname{conv}(\mathcal{X})$. 
By Definition~1 in~\citet{peters_convex_1986}, a function $f$ defined on a possibly nonconvex set $T$ is convex if, for any finite collection $\{x^j\}_{j=1}^n \subset T$ and any convex combination $\sum_{j=1}^n p_j x^j \in T$, the inequality
\begin{equation}
    \sum_{j=1}^n p_j\, f(x^j) \ge f\left( \sum_{j=1}^n p_j x^j \right)
\end{equation}
holds. Because every such convex combination of points in $\mathcal{X}$ lies in $\operatorname{conv}(\mathcal{X})$, where $\mathcal{Q}(x)$ is already known to be convex, the inequality is satisfied. 
Therefore, convexity follows directly for $\mathcal{X}$ once convexity over $\operatorname{conv}(\mathcal{X})$ is established. 
We highlight that, when $x$ is mixed-integer, the convexity of $\mathcal{Q}(x)$ trivially applies by the same reasoning.

\subsection{Case 3: mixed-integer first- and second-stage variables}

When both integer and continuous variables appear in the second stage, the recourse problem becomes the value function of a parametric mixed-integer linear program (MILP) in the parameter $x$. In contrast to the continuous settings of Cases 1 and 2, the value function $Q(x, \xi_s)$ of a parametric MILP is nonconvex in general. This is well established; see \citet{hassanzadeh_value_2014} for a detailed analysis of how the feasibility of integer assignments changes with the right-hand side parameter, producing nonconvex regions in the value function. An illustrative example is provided in~\ref{app:mixed_integer_dual}.

\subsubsection{A special case: quasi-convexity}
Even though convexity may fail, quasi-convexity may still arise in structured settings, and ICNNs often approximate it effectively in practice. To see this, consider the sublevel set
\begin{equation}
    L_{\alpha} = \{ x \mid Q(x, \xi_s) \le \alpha \},
\end{equation}
which are required to be convex for $Q(x, \xi_s)$ to be quasi-convex. In the mixed-integer setting, the recourse function for a fixed scenario can be written as
\begin{equation}
    Q(x, \xi_s) = \min_{y^{\mathrm{int}} \in \mathcal{Y}} \; Q_{y^{\mathrm{int}}}(x, \xi_s),
\end{equation}
where $y^{\mathrm{int}}$ denotes the second-stage integer variables and $\mathcal{Y}$ is the finite feasible set. 
For each fixed integer assignment $y^{\mathrm{int}}$, the remaining second-stage problem reduces to a continuous LP in the continuous recourse variables. The corresponding value function $Q_{y^{\mathrm{int}}}(x, \xi_s)$ is convex in $x$, as analysed in Cases 1 and 2.

Hence, the sublevel set
\begin{equation}
    L_{\alpha}
    = \bigcup_{y^{\mathrm{int}} \in \mathcal{Y}} 
      \left\{ x \mid Q_{y^{\mathrm{int}}}(x, \xi_s) \le \alpha \right\}
\end{equation}
is a union of convex sets.  
Thus, $Q(x, \xi_s)$ is quasi-convex whenever the union is convex, which occurs, for example, when:
\begin{enumerate}[label=(\roman*)]
    \item the sets $\left\{ x \mid Q_{y^{\mathrm{int}}}(x, \xi_s) \le \alpha \right\}$ are nested as $y^{\mathrm{int}}$ varies, or
    \item the feasible regions associated with different integer assignments produce LP subproblems whose sublevel sets overlap in a way that yields a convex union, or
    \item integer patterns activate in a monotone or threshold-dependent manner with respect to $x$, so that switching between convex regimes does not introduce nonconvexity in the union.
\end{enumerate}

Scenario-wise quasi-convexity does not guarantee quasi-convexity of $\mathcal{Q}(x)$ since quasi-convexity is not preserved under weighted summation. However, when the scenario-wise functions share similar convex regimes or threshold structure, the expected recourse may exhibit approximate quasi-convexity or near convexity in practice. In such settings, ICNNs are well-suited to approximate such behaviour because the convexity of their output in $x$ captures the dominant convex shape of each regime while smoothing over the transitions between regimes.

\subsection{Practical implications for ICNN surrogates}
In summary, for a broad class of 2SP problems, the recourse function $\mathcal{Q}(x)$ inherits a piecewise convex structure, even accommodating relatively challenging instances with strategic integer decisions in the first stage and operational adjustments (continuous) or logical switches (binary) in recourse actions.
Such structural regularity makes approximation via ICNNs advantageous, as they inherently preserve convexity while learning surrogate models of $\mathcal{Q}(x)$, enabling the embedding of the ICNN surrogates into the original 2SP for efficient solving. 
In cases where convexity cannot be analytically guaranteed, a practical approach is first to test whether the recourse function exhibits approximate convexity, for example, by training a lightweight ICNN and evaluating whether its predictive accuracy is acceptable. A symmetric check for approximate concavity is also possible by negating the training labels and assessing whether an ICNN fits well under this transformation.

\section{Methodology}\label{sec4}

This section presents our ICNN-enhanced framework for 2SP that directly leverages the convexity properties of $\mathcal{Q}(x)$. The ICNN surrogate of the recourse function is embedded into the first-stage problem to obtain a reformulated 2SP without introducing integer variables.

\subsection{Base architecture}\label{sec4_1}

Building on the extensive form~\eqref{eq:ef}, Neur2SP constructs a surrogate model $\hat{\mathcal{Q}}(x)$ to explicitly model the expected second-stage cost $\mathcal{Q}(x)$ through two sequential transformations:
\begin{enumerate}[label=(\roman*)]
    \item \textit{Scenario Encoding}: consists of the $x$-independent processing of each scenario $\xi_s$ through a shared network $\Psi^1$ to transform each scenario into a latent feature vector, followed by mean-aggregation and final encoding via another network $\Psi^2$:
    \begin{equation}
        \xi_\lambda = \Psi^2\left(\frac{1}{|S|}\sum_{s=1}^{|S|} \Psi^1(\xi_s)\right),
    \end{equation}
    where $\xi_\lambda$ serves as a compact representation of the entire scenario ensemble while preserving its statistical features.
    \item \textit{Decision Mapping}: represents a conventional ReLU NN $\Phi$ that combines first-stage decisions $x$ with scenario encoding $\xi_\lambda$ to predict $\mathbb{E}_\xi[Q(x,\xi)]$:
    \begin{equation}
        \hat{\mathcal{Q}}(x) = \Phi(x, \xi_\lambda ),
    \end{equation}
    without requiring separate evaluations for each scenario.
\end{enumerate}

While effective for general 2SP, this architecture becomes suboptimal for convex problems due to its MIP-based surrogate embeddings. The embedding of ReLU NNs as MIP constraints within optimisation problems is detailed in~\ref{app:mip}.
In our approach, we retain Neur2SP's core architecture for scenario processing while being specialised for convex 2SP models. The key difference lies in $\Phi$'s architecture: $\Phi$ is formulated as an ICNN when $Q(x,\xi_s)$ is convex in $x$ (per Section~\ref{sec3}).

\subsection{Surrogate modelling}\label{sec4_2}

For convex 2SP problems, we define $\Phi$ as:
\begin{equation}
    \Phi_{\text{ICNN}}(x,\xi_\lambda) = z_K,  
\end{equation}
where layer activations $z_k$ follow:
\begin{equation}
    z_{k+1} = \sigma ( W_k z_k + S_k [x,\xi_\lambda] + b_k), \quad k = {0, \dots, K-1},
\label{eq:fcicnn_xi}
\end{equation}
with $z_0=[x,\xi_\lambda]$, $W_0 \equiv 0$ and architectural constraints: (i) $W_k \geq 0$; and (ii) $\sigma(\cdot)$ is convex and non-decreasing. This ensures $\Phi_{\text{ICNN}}$ is convex in $x$, matching the convexity of $\mathcal{Q}(x)$ while preserving Neur2SP's ability to process arbitrary scenario sets through the encoder ($\Psi^1$ and $\Psi^2$).

In essence, the ICNN surrogate learns to map $(x,\xi_\lambda)$ pairs to $\mathcal{Q}(x)$ values, with $\xi_\lambda$ acting as a ``context vector'' that encapsulates the probability distribution formed by the scenario set~$S$. 
As such, it enables the ICNN to learn a family of convex functions in $x$, with each specific realisation of $\xi_\lambda$ generating a distinct convex function that models $\mathcal{Q}(x)$ for its corresponding scenario distribution.

The training process follows a supervised paradigm: 
(i) generate training pairs $(x^{(i)},\mathbb{E}_\xi[Q(x^{(i)},\xi)])$ by solving EF instances, where $i$ indexes the training examples;
(ii) train $\Psi^1$, $\Psi^2$ and $\Phi_{\text{ICNN}}$ jointly to minimise the mean squared error (MSE):
\begin{equation}
    \min_\Theta \; \mathbb{E}\left[{\left(\Phi_\text{ICNN} \left(x, \Psi^2\left(\tfrac{1}{|S|}\sum_{s=1}^{|S|} \Psi^1(\xi_s)\right)\right) 
    - \mathbb{E}_\xi[Q(x,\xi)]\right)}^2\right],
\end{equation}
where $\Theta = \{W_k, S_k, b_k\} \cup \theta_{\Psi^1} \cup \theta_{\Psi^2}$ represents the complete set of trainable parameters of $\Phi_{\text{ICNN}}$, $\Psi^1$, and $\Psi^2$.
Importantly, once trained, the model generalises to arbitrary second-stage scenario sets over the same problem structure without requiring retraining, eliminating scenario set-specific training overhead.

\subsection{Surrogate embedding}\label{sec4_3}

The ICNN architecture enables an exact convex optimisation formulation to be used instead of Neur2SP's MIP formulation. Building on the exact inference equivalence (Section~\ref{sec2_3}), we reformulate the 2SP as:
\begin{mini!}
    {x,z_k}{c^\intercal x + z_K}{\label{eq:exact-lp}}{}
    \addConstraint{z_{k+1}}{ \geq W_k z_k + S_k [x,\xi_\lambda] + b_k, \quad}{k=0, \ldots, K-1}
    \addConstraint{z_k} {\geq 0, \quad}{k=1,\ldots,K}
    \addConstraint{x \in \mathcal{X}.}
\end{mini!}
At deployment, $\xi_\lambda$ is precomputed offline from a given scenario set, ensuring the problem remains convex in $x$. 

This alternative embedding offers important methodological advantages. Compared to an MIP-based formulation, our approach significantly reduces the scale of mathematical programming formulation and the computational complexity. The conventional ReLU encoding necessitates introducing $N$ binary variables, $N$ continuous slack variables, and $N$ activation value variables (where $N$ is the neuron count). These auxiliary binary variables introduce combinatorial complexity that scales exponentially with network depth and width. In contrast, by leveraging ICNN's architectural properties, our formulation eliminates both these auxiliary integer variables and continuous slack variables entirely, yielding a more compact representation, which in turn benefits computational efficiency.

\section{Experiments}\label{sec5}

In this section, we report the results of our proposed ICNN-enhanced method across the three problem settings in terms of training and solving performance.

\subsection{Experimental setup}\label{sec5_1}

All experiments were conducted on the Triton computing cluster maintained by the Aalto Science-IT project, utilising an Intel Xeon Gold 6148 CPU and an NVIDIA Tesla H100 GPU with 64 GB of RAM. Our code extends Neur2SP’s codebase~\citep{khalil_neur2sp_2024} by augmenting it to include all necessary ICNN-related functionality. 
The implementation uses Python 3.11.5 with Gurobi 11.0.2~\citep{gurobi_optimization_llc_gurobi_2024} as the optimisation solver, and PyTorch 2.0.1~\citep{paszke_pytorch_2019} and scikit-learn 1.3.0~\citep{pedregosa_scikit_2011} for NN training. Full reproducibility materials are provided at~\url{https://github.com/Lycle/ICNN2SP}.

We evaluate our method on three canonical 2SP problems that are also used as benchmark problems in Neur2SP: the capacitated facility location problem (CFLP)~\citep{cornuejols_comparison_1991}, the stochastic server location problem (SSLP)~\citep{ntaimo_million-variable_2005}, and the investment problem (INVP)~\citep{schultz_solving_1998}. Problem instances are denoted by type followed by structural parameters, with the final parameter indicating scenario count. Detailed descriptions are provided in~\ref{app:spp}. 
The EF serves as the ground-truth baseline (under a 3-hour time limit). We compare both ICNN-enhanced 2SP and Neur2SP against EF to contextualise performance.

To ensure direct comparability with Neur2SP, we replicate their instance generation protocol, including uncertainty distributions and scenario sampling methods. We generate 5,000 first-stage solutions by uniformly sampling the feasible region. For each solution, the expected recourse cost is approximated by solving its second-stage problem for a random subset of scenarios, with the number of scenarios per solution randomly selected up to 100, which balances computational efficiency with the diversity of scenarios. Further details on data generation are provided in~\ref{app:dg}.
For both methods, we conduct a random search over 100 hyperparameter configurations for model selection, detailed in~\ref{app:hyper_rs}.
To mitigate solver-induced variability and isolate formulation-level effects, Gurobi’s parameters are fixed across all runs, as detailed in~\ref{app:paras}, following a similar practice as in~\citep{kronqvist_p-split_2025}.

\subsection{Results and discussion}\label{sec5_2}

Tables~\ref{tab:train_times} to~\ref{tab:obj_INVP} report our experimental results. 
We first examine training behaviour in terms of training time and validation accuracy (Table~\ref{tab:train_times}), then compare solving performance against EF using gap statistics and solving times (Tables~\ref{tab:res_CFLP}--\ref{tab:res_INVP}).
To complement these relative performance indicators, Tables~\ref{tab:obj_CFLP}--\ref{tab:obj_INVP} list the underlying objective values for all methods. 
In addition, we include results for the SSLP instances from the Stochastic Integer Programming Test Problem Library (SIPLIB)~\citep{ahmed_siplib_2015} in~\ref{app:siplib}. We also analyse the variable counts in~\ref{app:var} to compare formulation size.

\subsubsection{Training performance}

Both methods share the same structural design for Scenario Encoding ($\Psi^1$ and $\Psi^2$) while differing in the architectural choice of Decision Mapping ($\Phi$): our method uses ICNNs while Neur2SP employs standard ReLU NNs; this distinction constitutes the primary source of training performance differences.
As shown in Table~\ref{tab:train_times}, both methods exhibit training times within the same order of magnitude, with ICNNs requiring only marginally longer training due to their convexity constraints. 
This modest overhead is counterbalanced by their comparable or superior validation accuracy, with lower mean absolute error (MAE) achieved in 7 out of 10 problem instances.
These empirical findings are consistent with the convexity properties analysed in Section~\ref{sec3}.
In particular, ICNNs are architecturally constrained to represent convex functions, whereas standard ReLU networks used in Neur2SP can approximate arbitrary functions. 
The fact that ICNNs achieve validation accuracy comparable to or better than unconstrained networks suggests that the recourse functions in these benchmark problems can be well approximated by convex functions and that ICNNs can effectively capture their underlying structure, providing empirical evidence of a convex or near-convex recourse behaviour.
These results demonstrate that the convexity constraints do not significantly compromise the expressive power of the NN while presenting mathematical properties that enable more efficient downstream optimisation.

\begin{table}[ht]\centering
\footnotesize
\caption{Training times and MAE on validation set when applying ICNN-enhanced 2SP (ICNN) and Neur2SP (NN). Samples were split 80\%-20\% for train-validation, with the best model selected on validation set within the given epochs.}
\vspace{4pt}
\begin{tabular}{lrrrr}
\toprule
& \multicolumn{2}{c}{Training time (s)} & \multicolumn{2}{c}{Validation MAE} \\
\cmidrule(lr){2-3}\cmidrule(lr){4-5}
Problem  & {ICNN} & {NN} & {ICNN} & {NN}  \\
\midrule
CFLP\_10\_10 &  359.41 & 279.96 & 248.441 & 257.213\\
CFLP\_25\_25 &   288.77&  275.49 & 313.604 & 341.477 \\
CFLP\_50\_50 &   290.76&   270.52 & 558.073 & 569.384 \\
\midrule
SSLP\_5\_25 & 295.31 &   270.52 & 16.862 &  16.594 \\
SSLP\_15\_45 & 295.36 &    275.74 & 27.294 &  29.062\\
SSLP\_10\_50 & 292.08 &   273.90 & 17.372 & 14.105 \\
\midrule
INVP\_B\_E & 312.51 &    283.28 & 2.050   & 2.179\\
INVP\_B\_H & 312.04 & 285.02  & 1.902  & 2.031\\
INVP\_I\_E & 362.10 &     336.93 & 2.480 &  2.465 \\
INVP\_I\_H &  312.27 &     284.45 & 2.159 & 2.245\\
\bottomrule
\end{tabular}
\label{tab:train_times}
\end{table}

\subsubsection{Optimisation results}
Tables~\ref{tab:res_CFLP}--\ref{tab:res_INVP} reveal consistent performance advantages for ICNN-enhanced 2SP across problem classes. The reported standard deviations for solving times are generally small relative to the mean values, indicating limited variability across instances.
For CFLP, our method demonstrates superior scalability, achieving 1.3$\times$ speedups for small instances and up to 100$\times$ faster solving times for more complex problem class (CFLP\_50\_50) when compared to Neur2SP, while maintaining comparable or better solution quality. 
For SSLP, our approach achieves 2$\times$ speedups with identical optimality gaps for SSLP\_5\_25, a substantial gap reduction with 1.1--1.2$\times$ speedups for SSLP\_15\_45, and 3$\times$ faster solution times with matching solution quality for SSLP\_10\_50 instances.
For INVP, while absolute solving time differences are marginal ($\leq$ 0.03 s), ICNN-enhanced 2SP consistently delivers superior solution quality across most instances, often with considerable optimality gap improvements over Neur2SP.
As the scenario count increases, EF requires considerably more time to find solutions of comparable quality. For some of the largest scenario sets, EF fails to obtain solutions of the same quality, whereas our method obtains better quality solutions within seconds.

It is also worth noting that the absolute runtime differences between the Neur2SP results reported here and those in the original Neur2SP paper are primarily attributable to variations in hardware, software environments, solver parameters, and the model selection through hyperparameter search. 
A further factor is the formulation itself. The MIP-based embeddings of standard NNs introduce big-$M$ constraints and binary variables, which can lead to substantially higher variability and solver sensitivity. In contrast, the LP-based representation arising from the ICNNs avoids these factors and is expected to produce considerably more stable runtime behaviour across instances and solver configurations. 

\begin{table}[ht]\centering
\footnotesize
\caption{CFLP results (averaged over 11 instances per row; EF averages include only instances where EF found feasible solutions within 3h). Gap: relative objective difference (\%) of ICNN-enhanced 2SP (ICNN-EF) and Neur2SP (NN-EF) compared to EF; negative values indicate better solutions for minimisation. Solving time: runtime (s) reported as mean $\pm$ standard deviation; $\sim$10,800 indicates that the time limit was reached. EF time to match (s): average time for EF over instances when it achieved solution quality matching ICNN-enhanced 2SP (ICNN) and Neur2SP (NN); ``-'' indicates that EF did not match within the time limit; numbers in parentheses indicate the number of instances for which EF failed to match within the time limit. \textbf{Bold} indicates better performance.}
\vspace{4pt}
\resizebox{\textwidth}{!}{
\begin{tabular}{lrrrrrrr}
\toprule
  & \multicolumn{2}{c}{Gap (\%)} & \multicolumn{3}{c}{Solving time (s)} & \multicolumn{2}{c}{EF time to match (s) } \\
\cmidrule(lr){2-3}\cmidrule(lr){4-6}\cmidrule(lr){7-8}
Problem & ICNN-EF & NN-EF & ICNN & NN & EF & ICNN & NN\\
\midrule
CFLP\_10\_10\_100 & 0.57 & \textbf{0.51} & \textbf{0.29} $\pm$ 0.04 & 0.39 $\pm$ 0.07& 8,491.85 $\pm$ 3,840.29 & 3,801.86 \kern\digitwidth (0) & 5,281.23 \kern\digitwidth (0)\\
CFLP\_10\_10\_500 & 0.15 & \textbf{0.00} & \textbf{0.29} $\pm$ 0.03 & 0.39 $\pm$ 0.05& 10,800.05 $\pm$ 0.02 & 6,363.94 \kern\digitwidth (1) & 7,456.05 \kern\digitwidth (1)\\
CFLP\_10\_10\_1000 & -0.02 & \textbf{-0.27} & \textbf{0.29} $\pm$ 0.02 & 0.39 $\pm$ 0.02& 10,800.29 $\pm$ 0.53 & 5,946.09 \kern\digitwidth (3) & 6,331.67 \kern\digitwidth (3)\\
\midrule
CFLP\_25\_25\_100 & \textbf{2.19} & 3.37 & \textbf{0.39} $\pm$ 0.02 & 10.31 $\pm$ 0.12& 10,800.05 $\pm$ 0.01 & 327.97 \kern\digitwidth (2) & 235.64 \kern\digitwidth (0)\\
CFLP\_25\_25\_500 & \textbf{1.74} & 2.63 & \textbf{0.40} $\pm$ 0.02 & 10.37 $\pm$ 0.06& 10,800.12 $\pm$ 0.06 & 1,357.47 \kern\digitwidth (1) & 1,291.39 \kern\digitwidth (1)\\
CFLP\_25\_25\_1000 & \textbf{0.51} & 1.48 & \textbf{0.40} $\pm$ 0.02 & 10.07 $\pm$ 0.35& 10,800.33 $\pm$ 0.17 & 4,133.03 \kern\digitwidth (5) & 3,949.87 \kern\digitwidth (5)\\
\midrule
CFLP\_50\_50\_100 & \textbf{3.36} & 18.17 & \textbf{0.81} $\pm$ 0.02 & 84.07 $\pm$ 0.78 & 10,800.63 $\pm$ 1.49 & 1,432.75 \kern\digitwidth (0) & 317.61 \kern\digitwidth (0)\\
CFLP\_50\_50\_500 & \textbf{-9.82} & 3.06 & \textbf{0.84} $\pm$ 0.02 & 84.32 $\pm$ 0.12 & 10,800.16 $\pm$ 0.01 & -\kern\digitwidth (11) & 5,304.11\kern\digitwidth (10)\\
CFLP\_50\_50\_1000 & \textbf{-12.82} & -0.28 & \textbf{0.84} $\pm$ 0.03 & 84.39 $\pm$ 0.17 & 10,800.29 $\pm$ 0.01 & -\kern\digitwidth (11) & -\kern\digitwidth (11)\\
\bottomrule
\end{tabular}}
\label{tab:res_CFLP}
\end{table}

\begin{table}[ht]\centering
\footnotesize
\caption{SSLP results (averaged over 12 instances per row, including one from SIPLIB~\citep{ahmed_siplib_2015}). All settings and metric definitions follow Table~\ref{tab:res_CFLP}.}
\vspace{4pt}
\resizebox{\textwidth}{!}{
\begin{tabular}{lrrrrrrr}
\toprule
 & \multicolumn{2}{c}{Gap(\%)} & \multicolumn{3}{c}{Solving time (s)}  & \multicolumn{2}{c}{EF time to match (s)}\\
\cmidrule(lr){2-3}\cmidrule(lr){4-6}\cmidrule(lr){7-8}
Problem & ICNN-EF & NN-EF & ICNN & NN & EF & ICNN & NN\\
\midrule
SSLP\_5\_25\_50 & \textbf{0.12} & \textbf{0.12}  & \textbf{0.13} $\pm$ 0.02 & 0.25 $\pm$ 0.01& 1.58 $\pm$ 0.18 & 1.36 \kern\digitwidth (0) & 1.36 \kern\digitwidth (0) \\
SSLP\_5\_25\_100 & \textbf{0.02} & \textbf{0.02} & \textbf{0.12} $\pm$ 0.01 & 0.25 $\pm$ 0.01& 5.18 $\pm$ 0.34 & 5.01 \kern\digitwidth (0)  & 5.01 \kern\digitwidth (0)  \\
\midrule
SSLP\_15\_45\_5 & \textbf{2.65} & 18.09 & \textbf{0.16} $\pm$ 0.00 & 0.20 $\pm$ 0.01& 1.94 $\pm$ 2.26 & 1.15 \kern\digitwidth (0)& 0.11 \kern\digitwidth (0)\\
SSLP\_15\_45\_10 & \textbf{2.57} & 18.38 & \textbf{0.18} $\pm$ 0.01 & 0.20 $\pm$ 0.02& 32.74 $\pm$ 56.27 & 8.51 \kern\digitwidth (0)& 0.20 \kern\digitwidth (0)\\
SSLP\_15\_45\_15 & \textbf{2.11} & 18.72 & \textbf{0.20} $\pm$ 0.01 & 0.21 $\pm$ 0.02& 605.22 $\pm$ 1,103.58 & 538.02  \kern\digitwidth (0) & 0.40  \kern\digitwidth (0)\\
\midrule
SSLP\_10\_50\_50 & \textbf{0.00} & \textbf{0.00} & \textbf{0.27} $\pm$ 0.03 & 0.76 $\pm$ 0.01& 9,007.88 $\pm$ 3,950.65 & 12.28  \kern\digitwidth (0)& 12.28 \kern\digitwidth (0)\\
SSLP\_10\_50\_100 & \textbf{0.00} & \textbf{0.00} & \textbf{0.26} $\pm$ 0.08 & 0.76 $\pm$ 0.01& 10,800.04 $\pm$ 0.04 & 41.49  \kern\digitwidth (0)& 41.49 \kern\digitwidth (0)\\
SSLP\_10\_50\_500 & \textbf{0.00} & \textbf{0.00}  & \textbf{0.28} $\pm$ 0.07 & 0.77 $\pm$ 0.01& 10,800.21 $\pm$ 0.13 & 1,570.14 \kern\digitwidth (5)& 1,570.14 \kern\digitwidth (5)\\
SSLP\_10\_50\_1000 & \textbf{-2.22} & \textbf{-2.22} & \textbf{0.28} $\pm$ 0.11 & 0.77 $\pm$ 0.01& 10,803.79 $\pm$ 12.36 & -\kern\digitwidth (12) & -\kern\digitwidth (12)\\
SSLP\_10\_50\_2000 & \textbf{-102.68} & \textbf{-102.68} & \textbf{0.29} $\pm$ 0.10 & 0.78 $\pm$ 0.01& 10,800.14 $\pm$ 0.08 & -\kern\digitwidth (12)& -\kern\digitwidth (12)\\
\bottomrule
\end{tabular}}
\label{tab:res_SSLP}
\end{table}

\begin{table}[ht]\centering
\footnotesize
\caption{INVP results (one instance per row). All settings and metric definitions follow Table~\ref{tab:res_CFLP}.}
\vspace{4pt}
\begin{tabular}{lrrrrrrr}
\toprule
 & \multicolumn{2}{c}{Gap (\%)} & \multicolumn{3}{c}{Solving time (s)}  & \multicolumn{2}{c}{EF time to match (s)}\\
\cmidrule(lr){2-3}\cmidrule(lr){4-6}\cmidrule(lr){7-8}
Problem  & ICNN-EF & NN-EF & ICNN & NN & EF & ICNN & NN\\
\midrule
INVP\_B\_E\_4 & 11.06 & \textbf{10.80} & 0.06 & 0.15 & \textbf{0.00} & 0.00 & 0.00\\
INVP\_B\_E\_9 & 5.60  &  \textbf{4.31} & 0.07 & 0.10 & \textbf{0.01} & 0.00 & 0.00\\
INVP\_B\_E\_36 & \textbf{3.28} & 4.58 & 0.06 & 0.10 & \textbf{0.05} & 0.01 & 0.01\\
INVP\_B\_E\_121 & 3.07 & \textbf{2.85} & \textbf{0.06} & 0.09 & 1.32 & 0.01 & 0.01\\
INVP\_B\_E\_441 & \textbf{3.10} & 3.18 & \textbf{0.06} & 0.10 & 125.60 & 0.05 & 0.05\\
INVP\_B\_E\_1681 & 0.98 & \textbf{0.21} & \textbf{0.07} & 0.10 & 10,800.00 & 0.36 & 4,791.19\\
INVP\_B\_E\_10000 & -1.59  & \textbf{-1.64} & \textbf{0.07} & 0.11 & 10,808.33 & - & -\\
\midrule
INVP\_B\_H\_4 & \textbf{8.81} & \textbf{8.81} & 0.06 & 0.08 & \textbf{0.01} & 0.00 & 0.00\\
INVP\_B\_H\_9 & \textbf{5.04} & \textbf{5.04} & 0.07 & 0.09 & \textbf{0.01} & 0.01 & 0.01\\
INVP\_B\_H\_36 & \textbf{1.61} & \textbf{1.61} & \textbf{0.07} & 0.08 & 1.18 & 0.09 & 0.09\\
INVP\_B\_H\_121 & \textbf{1.77} & \textbf{1.77} & \textbf{0.06} & 0.08 & 60.09 & 0.06 & 0.06\\
INVP\_B\_H\_441 & \textbf{2.13} & \textbf{2.13} & \textbf{0.06} & 0.09 & 322.05 & 2.60 & 2.60\\
INVP\_B\_H\_1681 & \textbf{0.69} & \textbf{0.69} & \textbf{0.06} & 0.09 & 10,800.01 & 886.69 & 886.69\\
INVP\_B\_H\_10000 & \textbf{-2.72} & \textbf{-2.72} & \textbf{0.07} & 0.10 & 10,808.33 & - & -\\
\midrule
INVP\_I\_E\_4 & \textbf{4.78} & 17.32 & 0.05 & 0.08 & \textbf{0.00} & 0.00 & 0.00\\
INVP\_I\_E\_9 & \textbf{1.07} & 7.01 & 0.04 & 0.07 & \textbf{0.01} & 0.01 & 0.00\\
INVP\_I\_E\_36 & \textbf{0.37}  & 3.98 & 0.06 & 0.07 & \textbf{0.03} & 0.03 & 0.01\\
INVP\_I\_E\_121 & \textbf{0.17} & 2.88 & \textbf{0.06} & 0.08 & 0.39 & 0.20 & 0.01\\
INVP\_I\_E\_441 & \textbf{0.11} & 3.22 & \textbf{0.04} & 0.07 & 39.73 & 1.06 & 0.05\\
INVP\_I\_E\_1681 & \textbf{0.42}  & 3.39 & \textbf{0.05} & 0.07 & 10,800.00 & 4.18 & 0.64\\
INVP\_I\_E\_10000 & \textbf{0.22} & 3.50 & \textbf{0.05} & 0.08 & 10,800.09 & 670.21 & 465.22\\
\midrule
INVP\_I\_H\_4 & \textbf{13.78}  & \textbf{13.78} & 0.07 & 0.06 & \textbf{0.01} & 0.00 & 0.00\\
INVP\_I\_H\_9 & \textbf{5.85} & 9.12 & 0.05 & 0.07 & \textbf{0.01} & 0.01 & 0.01\\
INVP\_I\_H\_36 & \textbf{3.69} & 4.97 & \textbf{0.06} & 0.07 & 0.40 & 0.02 & 0.02\\
INVP\_I\_H\_121 & \textbf{3.75}  & 4.01 & \textbf{0.05} & 0.07 & 7.27 & 0.06 & 0.06\\
INVP\_I\_H\_441 & \textbf{1.27}  & 3.15 & \textbf{0.05} & 0.07 & 10,800.00 & 761.23 & 0.30\\
INVP\_I\_H\_1681 & \textbf{-0.41} & -0.34 & \textbf{0.05} & 0.07 & 10,800.01 & - & -\\
INVP\_I\_H\_10000 & 0.27  & \textbf{0.25} & \textbf{0.06} & 0.08 & 10,800.03 & 1,151.03 & 1,151.03\\
\bottomrule
\end{tabular}
\label{tab:res_INVP}
\end{table}

Tables~\ref{tab:obj_CFLP}--\ref{tab:obj_INVP} provide the corresponding objective values behind the gap statistics reported above.
For each problem class, the tables list the EF baseline objectives and the true and approximate objectives induced by the surrogate-based formulations. 
The true objective refers to the objective value obtained by taking the first-stage solution returned by the surrogate-based method and evaluating it in the original 2SP over the full scenario set. 
The approximate objective refers to the objective value obtained directly from the surrogate reformulation, where the expected second-stage cost is replaced by the learned surrogate model. 
The standard deviations are generally small relative to the corresponding means, suggesting that the observed performance patterns are consistent across instances. These results confirm that the ICNN-enhanced approach consistently achieves solution quality on par with conventional NN approaches, while retaining the computational advantages highlighted in Tables~\ref{tab:res_CFLP}--\ref{tab:res_INVP}.

\begin{table}[ht]\centering
\footnotesize
\caption{
Objectives for CFLP (averaged over 11 instances per row, reported as mean $\pm$ standard deviation). Baseline: EF objective at termination. True objective: cost of first-stage solutions from surrogate methods evaluated on original scenarios. Approximate objective: direct values from surrogate reformulations. ICNN: ICNN-enhanced 2SP, NN: Neur2SP. \textbf{Bold} indicates better solutions.}
\vspace{4pt}
\resizebox{\textwidth}{!}{
\begin{tabular}{lrrrrr}
\toprule
 &  & \multicolumn{2}{c}{True objective} & \multicolumn{2}{c}{Approximate objective} \\
\cmidrule(rr){3-4}\cmidrule(rr){5-6}
 Problem & Baseline & \multicolumn{1}{r}{ICNN} & \multicolumn{1}{r}{NN} & \multicolumn{1}{r}{ICNN}  & \multicolumn{1}{r}{NN} \  \\
\midrule
CFLP\_10\_10\_100 & 6,994.77 $\pm$ 150.27 & 7,034.16 $\pm$ 128.17 & \textbf{7,030.35} $\pm$ 156.59   & 7,172.48 $\pm$ 97.98 & 7,155.62 $\pm$ 126.00 \\
CFLP\_10\_10\_500 & 7,003.30 $\pm$ 47.79 & 7,013.85 $\pm$ 47.22 & \textbf{7,003.30} $\pm$ 47.79   & 7,179.69 $\pm$ 51.56 & 7,151.01 $\pm$ 70.34  \\
CFLP\_10\_10\_1000 & 7,013.16 $\pm$ 39.14 & 7,011.28 $\pm$ 53.62 & \textbf{6,994.31}  $\pm$ 17.11 & 7,158.10 $\pm$ 46.80 & 7,135.23 $\pm$ 42.17 \\
\midrule
CFLP\_25\_25\_100 & 11,854.23 $\pm$ 284.31 & \textbf{12,109.66} $\pm$ 112.18 & 12,252.05 $\pm$ 243.97  & 12,760.21 $\pm$ 0.00 & 11,396.54 $\pm$ 0.00  \\
CFLP\_25\_25\_500 & 11,856.48 $\pm$ 334.97 & \textbf{12,055.67} $\pm$ 37.60 & 12,160.15 $\pm$ 70.14   & 12,760.21 $\pm$ 0.00 & 11,396.54 $\pm$ 0.00 \\
CFLP\_25\_25\_1000 & 11,997.73 $\pm$ 466.37 & \textbf{12,047.44} $\pm$ 31.17 & 12,162.84 $\pm$ 60.64   & 12,760.21 $\pm$ 0.00 & 11,396.54 $\pm$ 0.00 \\
\midrule
CFLP\_50\_50\_100 & 25,573.72 $\pm$ 699.28 & \textbf{26,418.27} $\pm$ 318.30 & 30,205.61 $\pm$ 383.33  & 26,122.57 $\pm$ 0.00 & 23,556.44 $\pm$ 0.00 \\
CFLP\_50\_50\_500 & 29,352.39 $\pm$ 2,045.67 & \textbf{26,379.70} $\pm$ 107.34 & 30,172.44 $\pm$ 145.06   & 26,122.57 $\pm$ 0.00 & 23,556.44 $\pm$ 0.00  \\
CFLP\_50\_50\_1000 & 30,278.57 $\pm$ 271.92 & \textbf{26,395.75} $\pm$ 47.27 & 30,194.65 $\pm$ 99.91   & 26,122.57 $\pm$ 0.00 & 23,556.44 $\pm$ 0.00  \\
\bottomrule
\end{tabular}}
\label{tab:obj_CFLP}
\end{table}

\begin{table}[ht]\centering
\footnotesize
\caption{Objectives for SSLP (averaged over 12 instances per row). Same column definitions as in Table~\ref{tab:obj_CFLP}.}
\vspace{4pt}
\resizebox{\textwidth}{!}{
\begin{tabular}{lrrrrr}
\toprule
 & & \multicolumn{2}{c}{True objective} & \multicolumn{2}{c}{Approximate objective} \\
\cmidrule(lr){3-4} \cmidrule(lr){5-6}
 Problem & Baseline & ICNN & NN & \qquad ICNN  & NN \  \\
\midrule
SSLP\_5\_25\_50 & -125.36 $\pm$ 6.15  & \textbf{-125.22} $\pm$ 6.15 & \textbf{-125.22} $\pm$ 6.15   & -119.09 $\pm$ 0.00 & -124.47 $\pm$ 0.00  \\
SSLP\_5\_25\_100 & -120.94  $\pm$ 5.19  & \textbf{-120.91} $\pm$ 5.19 & \textbf{-120.91} $\pm$ 5.19   & -119.09 $\pm$ 0.00 & -124.47 $\pm$ 0.00  \\
\midrule
SSLP\_15\_45\_5 & -255.55 $\pm$ 19.14 & \textbf{-248.44} $\pm$ 15.63 & -208.49 $\pm$ 5.46 & -245.13 $\pm$ 0.00 & -227.99 $\pm$ 0.00  \\
SSLP\_15\_45\_10 & -257.41 $\pm$ 11.87  & \textbf{-250.62} $\pm$ 9.09 & -209.98 $\pm$ 2.77 & -245.13 $\pm$ 0.00 & -227.99 $\pm$ 0.00  \\
SSLP\_15\_45\_15 & -257.69 $\pm$ 6.79  & \textbf{-252.17} $\pm$ 4.24 & -209.08 $\pm$ 3.17  & -245.13 $\pm$ 0.00 & -237.45 $\pm$ 21.04  \\
\midrule
SSLP\_10\_50\_50 & -354.96 $\pm$ 7.41 & \textbf{-354.96} $\pm$ 7.41 & \textbf{-354.96} $\pm$ 7.41   & -357.70 $\pm$ 0.00 & -355.18 $\pm$ 0.00  \\
SSLP\_10\_50\_100 & -345.86 $\pm$ 5.91 & \textbf{-345.86} $\pm$ 5.91 & \textbf{-345.86} $\pm$ 5.91  & -357.70 $\pm$ 0.00 & -355.18 $\pm$ 0.00  \\
SSLP\_10\_50\_500 & -349.54 $\pm$ 2.62 & \textbf{-349.54} $\pm$ 2.62 & \textbf{-349.54} $\pm$ 2.62  & -357.70 $\pm$ 0.00 & -355.18 $\pm$ 0.00  \\
SSLP\_10\_50\_1000 & -342.64 $\pm$ 8.64  & \textbf{-350.07} $\pm$ 1.92 & \textbf{-350.07} $\pm$ 1.92 & -357.70 $\pm$ 0.00 & -355.18 $\pm$ 0.00  \\
SSLP\_10\_50\_2000 & -172.74 $\pm$ 2.39  & \textbf{-350.07} $\pm$ 1.80 & \textbf{-350.07} $\pm$ 1.80 & -357.70 $\pm$ 0.00 & -355.18 $\pm$ 0.00  \\
\bottomrule
\end{tabular}}
\label{tab:obj_SSLP}
\end{table}

\begin{table}[ht]\centering
\footnotesize
\caption{Objectives for INVP (one instance per row). Same column definitions as in Table~\ref{tab:obj_CFLP}.}
\vspace{4pt}
\begin{tabular}{lrrrrr}
\toprule
& & \multicolumn{2}{c}{True objective} & \multicolumn{2}{c}{Approximate objective} \\
\cmidrule(lr){3-4}\cmidrule(lr){5-6}
 Problem & Baseline & ICNN & NN & \qquad ICNN  & NN \  \\
\midrule
ip\_B\_E\_4 & -57.00 & -50.70 & \textbf{-50.84}   & -59.15 & -58.58 \\
ip\_B\_E\_9 & -59.33 & -56.01 & \textbf{-56.78}  & -59.63 & -59.27 \\
ip\_B\_E\_36 & -61.22 & \textbf{-59.21} & -58.42   & -59.69 & -59.30  \\
ip\_B\_E\_121 & -62.29 & -60.38 & \textbf{-60.51}   & -59.62 & -59.62  \\
ip\_B\_E\_441 & -61.32 & \textbf{-59.42} & -59.37   & -59.56 & -59.70  \\
ip\_B\_E\_1681 & -60.52 & -59.93 & \textbf{-60.39}    & -59.57 & -60.00  \\
ip\_B\_E\_10000 & -59.09 & -60.02 & \textbf{-60.06}   & -60.27 & -60.03  \\
\midrule
ip\_B\_H\_4 & -56.75 & \textbf{-51.75} & \textbf{-51.75}   & -55.70 & -59.17  \\
ip\_B\_H\_9 & -59.56 & \textbf{-56.56} & \textbf{-56.56}    & -57.22 & -59.11  \\
ip\_B\_H\_36 & -60.28 & \textbf{-59.31} & \textbf{-59.31}   & -59.38 & -59.58  \\
ip\_B\_H\_121 & -61.01 & \textbf{-59.93} & \textbf{-59.93}   & -59.72 & -59.96  \\
ip\_B\_H\_441 & -61.44 & \textbf{-60.14} & \textbf{-60.14}   & -60.50 & -60.07  \\
ip\_B\_H\_1681 & -60.89 & \textbf{-60.47} & \textbf{-60.47}   & -61.52 & -60.09  \\
ip\_B\_H\_10000 & -58.93 & \textbf{-60.53} & \textbf{-60.53}   & -61.07 & -60.29  \\
\midrule
ip\_I\_E\_4 & -63.50 & \textbf{-60.47} & -52.50   & -65.31 & -76.35  \\
ip\_I\_E\_9 & -66.56 & \textbf{-65.84} & -61.89    & -66.19 & -76.46  \\
ip\_I\_E\_36 & -69.86 & \textbf{-69.60} & -67.08   & -66.79 & -71.99  \\
ip\_I\_E\_121 & -71.12 & \textbf{-70.99} & -69.07   & -67.24 & -75.88  \\
ip\_I\_E\_441 & -69.64 & \textbf{-69.56} & -67.39   & -67.34 & -75.71 \\
ip\_I\_E\_1681 & -68.85 & \textbf{-68.56} & -66.52   & -67.48 & -75.82  \\
ip\_I\_E\_10000 & -68.05 & \textbf{-67.90} & -65.67   & -67.57 & -75.79  \\
\midrule
ip\_I\_H\_4 & -63.50 & \textbf{-54.75} & \textbf{-54.75}   & -65.10 & -65.33  \\
ip\_I\_H\_9 & -65.78 & \textbf{-61.93} & -59.78   & -63.71 & -65.48  \\
ip\_I\_H\_36 & -67.11 & \textbf{-64.63} & -63.78   & -65.57 & -66.44  \\
ip\_I\_H\_121 & -67.75 & \textbf{-65.21} & -65.03    & -66.99 & -66.76  \\
ip\_I\_H\_441 & -67.24 & \textbf{-66.38} & -65.12    & -66.99 & -66.76  \\
ip\_I\_H\_1681 & -65.41 & \textbf{-65.68} & -65.63   & -66.99 & -66.73  \\
ip\_I\_H\_10000 & -65.82 & -65.65 & \textbf{-65.66}   & -66.98 & -67.01 \\
\bottomrule
\end{tabular}
\label{tab:obj_INVP}
\end{table}

Overall, ICNN-enhanced 2SP demonstrates computational advantages across various problem settings, with solution quality matching or exceeding Neur2SP in 80\% of cases. Instances with lower validation MAE tend to yield better solution quality in subsequent optimisation. While training incurs modest overhead (10\% longer on average), this upfront investment enables acceleration during optimisation through LP-driven formulation. Taken together, these results position ICNN-enhanced 2SP as a scalable, high-fidelity alternative for convex or convexifiable 2SP problems.

\clearpage

\section{Conclusions}\label{sec6}
This work addresses a central computational bottleneck in 2SP by demonstrating that, for a broad and practically relevant class of convex recourse problems, ICNNs can serve as accurate surrogates embedded exactly via LP, with their architectural convexity constraints strategically harnessed rather than treated as limitations.
Across diverse problem classes, ICNN-enhanced 2SP outperforms existing MIP-based methodologies, with its computational advantages becoming more pronounced as problem scale increases.
The computational efficiency makes our approach viable for time-critical decision-making applications requiring frequent resolving, such as real-time power dispatch and production scheduling, where solution time directly impacts operational feasibility, and MIP-based alternatives would be prohibitively expensive. 

An important limitation of our approach is its convexity requirement. As such, our method complements Neur2SP by excelling in convex or convexifiable settings, whereas ReLU-based MIP formulations remain the most suitable choice for non-convex recourse. On the other hand, while our experiments focus on linear and mixed-integer recourse, the theoretical underpinnings of ICNNs suggest broader applicability to nonlinear convex settings, a promising direction for future work.

Additionally, the field of 2SP suffers from a scarcity of standardised, challenging benchmarks. Existing instances often lack the scale and complexity of modern industrial problems. Validating ICNN-enhanced 2SP on emerging applications will be critical to assessing its practical value. 
Finally, the interplay between ICNNs and problem structure, particularly in partially convex regimes, warrants deeper exploration, as hybrid architectures could unlock solutions to previously intractable stochastic programs. By bridging convex optimisation and deep learning, this work takes important further steps towards scalable, trustworthy decision-making under uncertainty.

\section*{Acknowledgements}
We thank the editors and anonymous reviewers for their constructive comments and valuable suggestions. We acknowledge the computational resources provided by the Aalto Science-IT project. Yu Liu gratefully acknowledges the support from the China Scholarship Council and thanks Yuanxiang Yang for insightful discussions.

\section*{Data availability}
All data and implementations are available in a repository referenced in the paper.

\appendix

\section{Mixed-integer recourse: an illustrative example}
\label{app:mixed_integer_dual}

This appendix illustrates why mixed-integer second-stage variables generally prevent the existence of a single $x$-independent dual polyhedron for the recourse function, which in turn explains why convexity may fail even when each continuous subproblem is convex.

Consider a second-stage problem with one binary variable $y^{\mathrm{int}} \in \{0,1\}$ and one continuous variable $y^{\mathrm{cont}} \ge 0$:
\begin{mini!}
    {y^{\mathrm{int}},\, y^{\mathrm{cont}}}
    {q_1\, y^{\mathrm{int}} + q_2\, y^{\mathrm{cont}}}
    {}
    {}
    \addConstraint{y^{\mathrm{cont}}}{\ge h - T x - M y^{\mathrm{int}},}
    \addConstraint{y^{\mathrm{int}}}{\in \{0,1\},}
    \addConstraint{y^{\mathrm{cont}}}{\ge 0,}
\end{mini!}
where $M$ is a sufficiently large constant.

For a fixed integer assignment, the continuous subproblem is an LP with feasibility condition
\begin{equation}
    \begin{cases}
        y^{\mathrm{cont}} \ge h - T x, & \text{if } y^{\mathrm{int}} = 0,\\[6pt]
        y^{\mathrm{cont}} \ge h - T x - M, & \text{if } y^{\mathrm{int}} = 1.
    \end{cases}
\end{equation}
Hence, the assignment $y^{\mathrm{int}} = 1$ is feasible only when $h - T x - M \le 0$, showing that the feasibility of integer patterns depends on $x$.

The induced dual regions also differ. Let $\pi_0$ and $\pi_1$ denote the dual variables associated with the constraints in the two cases. The dual feasible region associated with each integer assignment can be written as
\begin{equation}
    \begin{cases}
        D_0 = \{ \pi_0 \ge 0 \}, & \text{if } y^{\mathrm{int}} = 0,\\[6pt]
        D_1 = \{ \pi_1 \ge 0 \}, & \text{if } y^{\mathrm{int}} = 1.
    \end{cases}
\end{equation}

The recourse value is generated by $D_0$ when $h - T x > M$ and by $D_1$ when $h - T x \le M$. Because the active dual region depends on $x$, no single dual polyhedron $D$ satisfies
\begin{equation}
    Q(x, \xi_s) = \max_{\pi \in D} \, \bigl\{ -\pi^{\intercal} T x + \pi^{\intercal} h \bigr\}
    \qquad
\end{equation}
for all $x$. This is in contrast with the continuous recourse case, where an $x$-independent dual region exists and the recourse function is a pointwise maximum of affine functions. With mixed-integer recourse, this fixed dual representation no longer exists, so the max-of-affine representation breaks down and convexity in $x$ may fail.

\section{MIP reformulations}
\label{app:mip}
This appendix details the MIP formulation required to embed conventional ReLU NNs into optimisation problems, highlighting the additional complexity that our ICNN-based approach avoids.

Consider a conventional ReLU NN with $L+1$ layers (numbered from $0$ to $L$) to build the surrogate model for an underlying function $f$, which can be approximated and effectively replaced using this ReLU NN surrogate $\hat f$, constructed based on dataset $\mathcal{D}$. For the input layer, we have $y^0 = x$. For the output layer we have $y^L = \hat{f}_\mathcal{D}(x) $. The activations of neurons $i = 1, \dots, N_l$ in the hidden layer $l = 1, \dots, L-1$ are calculated as
\begin{equation}
    y^l_i = \text{ReLU}(w_i^{l^\intercal} y_{l-1} + b^l_i) = \max\,(0,w_i^{l^\intercal} y_{l-1} + b^l_i),
\end{equation}
where $w_i^l$ and $b^l_i$ denote the weight and bias of the corresponding neuron, respectively. Suppose we could obtain the lower bounds $L_i^l < 0$ and upper bounds $U_i^l > 0$ such that
\begin{align}
L_i^l \leq w_i^{l^\intercal} y_{l-1} + b^l_i \leq U_i^l.
\end{align}

Following~\citet{fischetti_deep_2018}, the ReLU operator can be encoded into mixed-integer linear constraints by introducing slack variables $s^l_i$ and binary variables $z^l_i$ for $i = 1, \dots, N_l,\ l = 1, \dots, L-1$, which is expressed as:
\begin{align}
    & y_{i}^{l}-s_{i}^{l} = w_{i}^{l^{\intercal}}y^{l-1}+b_{i}^{l}, \\
    & 0\leq y_{i}^{l}\leq U_{i}^{l}z_{i}^{l}, \\
    & 0\leq s_{i}^{l}  \leq -L_{i}^{l}\left( 1-z_{i}^{l} \right) ,\\
    & z_{i}^{l}  \in \left\{ 0,1 \right\}.
\end{align}

Combining these constraints with input and output layer bounds 
\begin{align}
    & L^0 \leq y^0 \leq U^0, \\
    & L^L \leq y^L =  w^{L^{\intercal}}y^{L-1}+b^{L} \leq U^L,
\end{align}
we obtain the exact MIP model to embed the ReLU NN surrogate $ y = \hat{f}_\mathcal{D}(x)$ into a larger optimisation problem.
This formulation necessarily introduces $N$ binary variables, $N$ continuous slack variables, and $N$ activation value variables, where $N=\sum_{l=1}^{L-1} N_l$ represents the total number of neurons across all hidden layers.
The inclusion of these binary variables fundamentally alters the computational complexity class of the resulting optimisation problem, leading to exponential worst-case solution time as network size increases.

\section{Stochastic programming problems}
\label{app:spp}

This appendix provides brief descriptions of the three benchmark problems employed in our experimental evaluation. These are the same problems utilised in \citet{dumouchelle_neur2sp_2022}, facilitating direct comparison between our ICNN-enhanced approach and the Neur2SP methodology. The pooling problem (PP) considered in~\citet{dumouchelle_neur2sp_2022} is not included in this work because it contains bilinear blending constraints that induce a nonlinear and generally nonconvex recourse function. Since the proposed approach relies on ICNN surrogates that enforce convexity with respect to the first-stage decisions to enable an exact LP embedding, the pooling problem falls outside our method's modelling assumptions.

\paragraph{Capacitated facility location problem (CFLP)}
The CFLP addresses the challenge of selecting optimal facility locations to serve customer demand at minimal cost. In its stochastic formulation, facility decisions precede knowledge of actual demand patterns. Our implementation adapts the case from \citet{cornuejols_comparison_1991} to the stochastic setting by preserving the original cost structures and capacity parameters while stochastically sampling demand values. To guarantee relatively complete recourse, we incorporate penalty variables for potential demand shortfalls. Throughout our experimental analysis, we denote a CFLP instance with $n$ facilities, $m$ customers, and $s$ scenarios as CFLP\_$n$\_$m$\_$s$.

\paragraph{Stochastic server location problem (SSLP)}
The SSLP focuses on optimal server deployment and client assignment under uncertain service requests. Binary first-stage variables determine server placement, while second-stage variables manage client-server assignments after observing which clients require service. Our experimental protocol utilises the standard benchmark instances from SIPLIB \citet{ahmed_siplib_2015} as documented in \citet{ntaimo_million-variable_2005}. For consistent reference, we denote an SSLP instance with $n$ potential server locations, $m$ possible clients, and $s$ scenarios as SSLP\_$n$\_$m$\_$s$. 

\paragraph{Investment problem (INVP)}
The INVP consists of continuous first-stage investment decisions followed by discrete second-stage profit-generating activities. We implement Example 7.3 from \citet{schultz_solving_1998}, with two first-stage variables bounded by $[0, 5]$ and four binary second-stage variables. Scenarios are derived from two random variables uniformly distributed over $[5,15]$. For consistency, we reformulate it as a minimisation problem. In our experiments, we denote instances as INVP\_$v$\_$t$\_$s$, where $v$ indicates second-stage variable type (B for binary, I for integer), $t$ indicates technology matrix type (E for identity, H for $[[2/3, 1/3], [1/3, 2/3]]$), and $s$ indicates the number of scenarios.

\section{Data generation}
\label{app:dg}

Table~\ref{tab:dg_times} reports computational time for generating the shared 5,000-sample dataset used to train both NN and ICNN surrogates for the expected recourse function. 
Under the dataset generation setting as described in Section~\ref{sec5_1}, for each of the 5,000 first-stage solutions, we randomly sample between 1 and 100 scenarios (averaging $\sim$50 scenarios per solution), resulting in approximately 5,000 $\times$ 50 = 250,000~second-stage optimisation problems to solve for each problem instance. 
Notably, this process is highly parallelisable and could be extended to utilise more computational resources when available.

\begin{table}[!ht]\centering
\footnotesize
\caption{Data generation time (5,000 samples per instance) utilising 16 parallel processes.}
\vspace{4pt}
\begin{tabular}{lrrr}
\toprule
 & \multicolumn{2}{c}{(IC)NN} \\
\cmidrule{2-3}
Problem & Time per sample (s)  & Total time (s)  \\
\midrule
CFLP\_10\_10  & 0.32 & 1,620.07 \\
CFLP\_25\_25  & 2.22 & 11,096.88  \\
CFLP\_50\_50  & 2.19 & 10,945.11 \\
\midrule
SSLP\_5\_25   & 0.84 & 4,182.60 \\
SSLP\_15\_45  & 2.49 & 12,438.02  \\
SSLP\_10\_50  & 2.33 & 11,649.44 \\
\midrule
INVP\_B\_E    & 1.77 & 8,855.89  \\
INVP\_B\_H    & 1.43 & 7,157.07 \\
INVP\_I\_E    & 1.48 & 7,420.66 \\
INVP\_I\_H    & 1.04 & 5,194.90  \\
\bottomrule
\end{tabular}
\label{tab:dg_times}
\end{table}

\section{Model selection}
\label{app:hyper_rs}

For model selection, we conducted a random search over 100 configurations for each problem instance following~\citet{dumouchelle_neur2sp_2022}. Both the ICNN-enhanced 2SP and Neur2SP methods were subjected to identical hyperparameter sampling from the distribution space outlined in Table~\ref{tab:rs}, with training proceeding for 2,000 epochs across all configurations. 

To ensure fair comparison, we constrained the architecture of both methods to utilise a single ReLU hidden layer with varying dimensionality within $\Phi$.
For the Scenario Encoding component, which precedes surrogate training as described in Section~\ref{sec4_1}, we implemented a three-layer architecture. Specifically, ``Embed hidden dimension 1'' and ``Embed hidden dimension 2'' correspond to the layers before scenario aggregation (within $\Psi^1$), while ``Embed hidden dimension 3'' represents the final hidden layer after aggregation (within $\Psi^2$). This configuration allows the model to effectively distil scenario information into the fixed vector representation $\xi_\lambda$. The best configurations obtained from the random search for each problem instance are reported in Tables~\ref{tab:icnn_best} and~\ref{tab:nn_best} for the ICNN-enhanced 2SP and Neur2SP models, respectively.

\begin{table}[!ht]
\centering
\footnotesize
\caption{NN hyperparameter sampling space. Values in \{\}: equal probability sampling; values in []: uniform distribution sampling.}
\vspace{4pt}
\begin{tabular}{lcc}
\toprule
Parameter                   &  Sampling range \\
\midrule
Batch size                  & $\{16, 32, 64, 128\}$  \\
Learning rate               & $[1e^{-5},1e^{-1}]$ \\
L1 weight penalty           & $[1e^{-5}, 1e^{-1}]$  \\
L2 weight penalty           & $[1e^{-5}, 1e^{-1}]$\\
Optimiser                   & \{Adam, Adagrad, RMSprop\}\\
Dropout                     & $[0, 0.5]$ \\     
ReLU hidden dimension       & $\{64, 128, 256, 512\}$  \\
Embed hidden dimension 1    & $\{64, 128, 256, 512\}$   \\
Embed hidden dimension 2    & $\{16, 32, 64, 128\}$  \\
Embed hidden dimension 3    & $\{8, 16, 32, 64\}$  \\
\bottomrule
\end{tabular}
\label{tab:rs}
\end{table}

\begin{table}[!ht]
\centering
\caption{Best configurations identified by random search for the ICNN-enhanced 2SP models.}
\label{tab:icnn_best}
\vspace{4pt}
\resizebox{\textwidth}{!}{
\begin{tabular}{lrrrrrrrrrr}
\toprule
{Parameter} & {CFLP\_10\_10} & {CFLP\_25\_25} & {CFLP\_50\_50} & {SSLP\_5\_25} & {SSLP\_15\_45} & {SSLP\_10\_50} & {INVP\_B\_E} & {INVP\_B\_H} & {INVP\_I\_E} & {INVP\_I\_H} \\
\midrule
    Batch size              & 128     & 128     & 128     & 128     & 16      & 16      & 32      & 32      & 32      & 32      \\    
    Learning rate           & 0.00436 & 0.08384 & 0.08384 & 0.08384 & 0.04383 & 0.02639 & 0.00768 & 0.00768 & 0.00768 & 0.00768 \\    
    L1 weight penalty       & 0.09067 & 0.03226 & 0.03226 & 0.03226 & 0.00976 & 0.00120& 0       & 0       & 0       & 0       \\    
    L2 weight penalty       & 0.02603 & 0.07454 & 0.07454 & 0.07454 & 0       & 0       & 0       & 0       & 0       & 0       \\    
    Optimizer               & Adam    & RMSprop & RMSprop & RMSprop & RMSprop & Adam    & RMSprop & RMSprop & RMSprop & RMSprop \\     
    Dropout                 & 0.04226 & 0.00238 & 0.00238 & 0.00238 & 0.04066 & 0.02918 & 0.07346 & 0.07346 & 0.07346 & 0.07346 \\    
    ReLU hidden dimension   & 512     & 512     & 512     & 512     & 256     & 512     & 256     & 256     & 256     & 256     \\  
    Embed hidden dimension  & 512     & 64      & 64      & 64      & 128     & 64      & 128     & 128     & 128     & 128     \\ 
    Embed dimension 1       & 64      & 128     & 128     & 128     & 16      & 16      & 64      & 64      & 64      & 64      \\ 
    Embed dimension 2       & 16      & 8       & 8       & 8       & 32      & 8       & 32      & 32      & 32      & 32      \\ 
\bottomrule
\end{tabular}}
\end{table}

\begin{table}[!ht]
\centering
\caption{Best configurations identified by random search for the Neur2SP models.}
\label{tab:nn_best}
\vspace{4pt}
\resizebox{\textwidth}{!}{
\begin{tabular}{lrrrrrrrrrr}
\toprule
{Parameter} & {CFLP\_10\_10} & {CFLP\_25\_25} & {CFLP\_50\_50} & {SSLP\_5\_25} & {SSLP\_15\_45} & {SSLP\_10\_50} & {INVP\_B\_E} & {INVP\_B\_H} & {INVP\_I\_E} & {INVP\_I\_H} \\
\midrule
    Batch size              & 128     & 128     & 128     & 128     & 128     & 128     & 128     & 128     & 128     & 128     \\    
    Learning rate           & 0.00436 & 0.08384 & 0.08384 & 0.08384 & 0.09621 & 0.08384 & 0.00433 & 0.00433 & 0.00433 & 0.00433 \\    
    L1 weight penalty       & 0.09067 & 0.03226 & 0.03226 & 0.03226 & 0.02965 & 0.03226 & 0.00500& 0.00500& 0.00500& 0.00500\\    
    L2 weight penalty       & 0.02603 & 0.07454 & 0.07454 & 0.07454 & 0.00040 & 0.07454 & 0.00841 & 0.00841 & 0.00841 & 0.00841 \\    
    Optimizer               & Adam    & RMSprop & RMSprop & RMSprop & Adam    & RMSprop & Adam    & Adam    & Adam    & Adam    \\     
    Dropout                 & 0.04226 & 0.00238 & 0.00238 & 0.00238 & 0.01053 & 0.00238 & 0.04056 & 0.04056 & 0.04056 & 0.04056 \\    
    ReLU hidden dimension   & 512     & 512     & 512     & 512     & 128     & 512     & 128     & 128     & 128     & 128     \\  
    Embed hidden dimension  & 512     & 64      & 64      & 64      & 256     & 64      & 512     & 512     & 512     & 512     \\ 
    Embed dimension 1       & 64      & 128     & 128     & 128     & 64      & 128     & 128     & 128     & 128     & 128     \\ 
    Embed dimension 2       & 16      & 8       & 8       & 8       & 64      & 8       & 8       & 8       & 8       & 8       \\ 
\bottomrule
\end{tabular}}
\end{table}

\section{Solver parameter settings}
\label{app:paras}
Gurobi’s automatic parameter selection introduces variability in results due to divergent solution strategies rather than formulation differences. To isolate methodological comparisons, we fix key Gurobi parameters (Table~\ref{tab:gurobi_paras}), standardising the solving environment across formulations. See~\citet{gurobi_parameter_2024} for details.

\begin{table}[!ht]\centering
\footnotesize
\caption{Solver parameter settings in the experiments.}
\vspace{4pt}
\begin{tabular}{lrl}
\toprule
Parameter & Value & Description \\
\midrule
\texttt{TimeLimit}  & 10800 & Limit total runtime (3 h) \\
\texttt{Threads}    & 1 & Thread count \\
\texttt{Presolve}   & 0 & Disable presolve \\
\texttt{Heuristics} & 0 & Minimise time spent in MIP heuristics \\
\texttt{NodeMethod} & 1  & Use dual simplex for MIP node relaxations \\
\texttt{Cuts}       & 1 & Moderate cut generation \\
\texttt{MIPFocus}   & 2 & Focus on proving optimality \\
\bottomrule
\end{tabular}
\label{tab:gurobi_paras}
\end{table}

\section{SSLP SIPLIB results}
\label{app:siplib}

This section presents performance metrics for the SSLP instances from SIPLIB~\citep{ahmed_siplib_2015}. Table~\ref{tab:res_siplib} reports optimality gaps and solving times across different problem configurations.

The results reveal that ICNN-enhanced 2SP achieves optimal or near-optimal solutions (0--1.63\% gaps) across all instances. For medium-sized SSLP\_15\_45 instances, ICNN-enhanced 2SP significantly outperforms Neur2SP in solution quality while maintaining comparable solving times. This notable performance divergence suggests that the ICNN architecture better captures the underlying structure of these specific problem instances. 
Computational efficiency advantages of ICNN-enhanced 2SP increase with problem complexity. For smaller instances (SSLP\_5\_25), ICNN shows approximately 2$\times$ speedup over Neur2SP. For the largest instances (SSLP\_10\_50), both surrogate methods vastly outperform EF, while ICNN-enhanced 2SP consistently requires approximately one-third the solving time of Neur2SP for these challenging instances.

\begin{table}[!ht]\centering
\footnotesize
\caption{
Performance comparison for SSLP SIPLIB instances. ICNN: ICNN-enhanced 2SP, NN: Neur2SP, EF: Extensive Form (baseline). \textbf{Bold} indicates superior performance between ICNN and NN (lower is better). Optimal reference values are from \citet{ahmed_siplib_2015}.
}
\vspace{4pt}
\begin{tabular}{lrrrrrr}
\toprule
 & \multicolumn{3}{c}{Gap to optimum (\%)} & \multicolumn{3}{c}{Solving time (s)} \\
\cmidrule(lr){2-4}\cmidrule(lr){5-7}
Problem & ICNN & NN & EF  & ICNN & NN & EF\  \\
\midrule
SSLP\_5\_25\_50 & \textbf{0.00} & \textbf{0.00} & 0.00 & \textbf{0.12} & 0.24 & 1.46 \\
SSLP\_5\_25\_100 & \textbf{0.00} & \textbf{0.00} & 0.00 & \textbf{0.12} & 0.26 & 5.44 \\
\midrule
SSLP\_15\_45\_5 & \textbf{0.53} & 18.90 & 0.00 & \textbf{0.16} & 0.19 & 3.36 \\
SSLP\_15\_45\_10 & \textbf{1.54} & 18.23 & 0.00 & \textbf{0.16} & 0.20 & 1.57 \\
SSLP\_15\_45\_15 & \textbf{1.63} & 16.51 & 0.00   & 0.28 & \textbf{0.20} & 7.14 \\
\midrule
SSLP\_10\_50\_50 & \textbf{0.00} & \textbf{0.00} & 0.00 & \textbf{0.27} & 0.77 & 10,800.02 \\
SSLP\_10\_50\_100 & \textbf{0.00} & \textbf{0.00} & 0.00 & \textbf{0.26} & 0.76 & 10,800.16 \\
SSLP\_10\_50\_500 & \textbf{0.00} & \textbf{0.00} & 0.00 & \textbf{0.24} & 0.76 & 10,800.59 \\
SSLP\_10\_50\_1000 & \textbf{0.00} & \textbf{0.00} &  4.77 & \textbf{0.27} & 0.75 & 10,800.08 \\
SSLP\_10\_50\_2000 & \textbf{0.00} & \textbf{0.00} & 51.24 & \textbf{0.27} & 0.76 & 10,800.10 \\
\bottomrule
\end{tabular}
\label{tab:res_siplib}
\end{table}

\section{Formulation size}
\label{app:var}

We compare variable counts of the surrogate embedding reformulation between ICNN-enhanced 2SP (our method) and Neur2SP (MIP-based method) in Table~\ref{tab:var_nums}. As detailed in Section~\ref{sec4}, ICNN’s convex structure avoids introducing new integer variables and halves the auxiliary continuous variables required by Neur2SP’s ReLU activation encodings. 

\begin{table}[!ht]\centering
\footnotesize
\caption{Variable counts comparison between ICNN-enhanced 2SP and Neur2SP reformulations. \textit{Italic} entries indicate instances with architecture parity after hyperparameter tuning.}
\vspace{4pt}
\begin{tabular}{lrrrr}
\toprule
 & \multicolumn{2}{c}{\# Continuous} & \multicolumn{2}{c}{\# Integer (Binary)}  \\
\cmidrule(lr){2-3}\cmidrule(lr){4-5}
Problem & {ICNN} & {NN} & {ICNN} & {NN} \\
\midrule
\textit{CFLP\_10\_10}&  513&  1025&    10&  522\\
\textit{CFLP\_25\_25} &   513&  1025&   25&  537\\
\textit{CFLP\_50\_50} &   513&  1025&   50&  562\\
\midrule
\textit{SSLP\_5\_25} &   513& 1025&   5&    517\\
SSLP\_15\_45 &   257&  257&   15&  143\\
\textit{SSLP\_10\_50} &   513&  1025&   10&  522\\
\midrule
INVP\_B\_E &  259 &   259&  0 &      128\\
INVP\_B\_H &  259 &   259& 0  &      128\\
INVP\_I\_E &  259 &   259& 0  &      128\\
INVP\_I\_H &  259 &   259& 0  &     128\\
\bottomrule
\end{tabular}
\label{tab:var_nums}
\end{table}

For integer variables: ICNN uses only original first-stage integers (e.g., in CFLP/SSLP; none in INVP); Neur2SP adds auxiliary integers to encode ReLU activations.
For continuous variables: ICNN requires 50\% fewer auxiliary variables than Neur2SP with equivalent network architectures. This reduction stems from ICNN's LP-representable structure, which requires only one additional variable per neuron, whereas ReLU networks in Neur2SP require two (one for the activation output and one slack variable).
Network architecture impacts reformulation complexity. \textit{Italic} entries mark instances where both methods share identical architectures (same hidden layer dimension) after hyperparameter tuning. In these cases, ICNN's reformulation efficiency is most pronounced, with variable count reductions directly attributable to methodological differences rather than model size variations.

{
\small
\bibliographystyle{elsarticle-num-names}
\bibliography{references}
}

\end{document}